\documentclass[english,11pt]{article}

\usepackage[latin1]{inputenc}
\usepackage{amsmath,amsthm,amssymb}
\usepackage{graphicx}
\usepackage{color}

\def\proof{{\it Proof.}\ }

\def\eq#1{(\ref{#1})}

\def\neweq#1{\begin{equation}\label{#1}}
\def\endeq{\end{equation}}

\def\la{\lambda}

\def\RR{{\mathbb R} }

\def\di{\displaystyle}
\def\ri{\rightarrow}

\def\ii{\^\i }
\newtheorem{theorem}{Theorem}[section]
\newtheorem{lem}{Lemma}[section]

\title{\sc Ground state solutions for the singular Lane-Emden-Fowler  
equation 
with sublinear convection term}
\author{Marius GHERGU$^a$ and
Vicen\c tiu R\u ADULESCU$^b$\\
\small $^a$ Institute of Mathematics ``Simion Stoilow" of the Romanian 
Academy,\\ 
\small 21, Calea Grivitei Street, 010702 Bucharest, Sector 1, Romania\\
\small E-mail: {\tt marius.ghergu@imar.ro}\\
\small $^b$ Department of Mathematics, University of Craiova,
200585 Craiova, Romania\\ \small E-mail: {\tt
vicentiu.radulescu@math.cnrs.fr}}
\date{}

\begin{document}
\baselineskip16pt \maketitle
\renewcommand{\theequation}{\arabic{section}.\arabic{equation}}
\catcode`@=11 \@addtoreset{equation}{section} \catcode`@=12

\bigskip
\begin{center}{This paper is dedicated to Professor William F. Ames on 
his 80th birthday}
\end{center}

\vspace{2cm}
\noindent {\bf Correspondence address:}\\
Vicen\c tiu R\u ADULESCU\\
Department of Mathematics\\ University of Craiova\\
200585 Craiova, Romania\\  E-mail: {\tt
vicentiu.radulescu@u-picardie.fr}
\newpage
\begin{abstract}
We are concerned with singular elliptic equations of the form
$-\Delta u= p(x)(g(u)+ f(u)+|\nabla u|^a)$ in $\RR^N$ ($N\geq 3$),
where $p$ is a positive weight and
$0<a<1$. Under the hypothesis that $f$ is a nondecreasing function with 
sublinear growth and $g$ is  decreasing and
unbounded around the origin, we establish the
existence of a ground state solution vanishing at infinity. Our 
arguments 
rely essentially on the maximum principle.
\end{abstract}

\section{Introduction and the main result}
The Lane-Emden-Fowler equation originated from earlier theories concerning gaseous dynamics in
astrophysics around the turn of the last century (see, e.g.,  Fowler \cite{fowler}). It also %%@
arises in the study of
fluid mechanics, relativistic mechanics, nuclear physics and in the study of chemical reaction
systems, one can see the survey article by Wong \cite{wong} for detailed background of the %%@
generalized
Emden-Fowler equation. The Lane-Emden-Fowler equation has been
studied by many authors using various methods and techniques. For example, critical point theory,
fixed point theory, topological degree theory and coincidence degree theory are widely used
to study the existence of the BVP for ordinary and partial differential equations (see, e.g., %%@
Agarwal and O'Regan \cite{ar}, Arcoya \cite{arcoya}, Hale and Mawhin \cite{hale}, Mawhin and %%@
Willem \cite{mw}, Nato and Tanaka \cite{nato}, Wang \cite{wang}, etc.).

We are concerned in this paper with the following singular
Lane-Emden-Fowler type problem 
\neweq{pp}
\left\{\begin{tabular}{ll}
$-\Delta u=p(x)(g(u)+f(u)+ |\nabla u|^a)$
\quad & $\mbox{\rm in}\ \RR^N,$\\
$u>0$ \quad & $\mbox{\rm in}\ \RR^N,$\\
$u(x)\ri 0$ \quad & $\mbox{\rm as}\ |x|\ri \infty,$
\end{tabular} \right.
\endeq
where $N\geq 3$, $0<a<1$, and $p:\RR^N\ri (0,\infty)$ is a H\"older 
continuous function of exponent $\gamma\in(0,1)$.
We assume that $g\in C^1(0,\infty)$ is a positive decreasing function
such that 

\medskip
\noindent $\di(g1)\qquad\lim_{t\ri 0^+}g(t)=+\infty.$

\medskip

Throughout this paper we suppose that
$f:[0,\infty)\rightarrow[0,\infty)$ is a
H\"{o}lder continuous function of exponent $0<\gamma<1$ which is
nondecreasing with respect to the second variable and
such that $f$ is positive on
$\overline\Omega\times (0,\infty).$
The analysis we
develop in this paper concerns the case where $f$ is sublinear, that 
is, 
\medskip

\noindent $\di(f1)\qquad$ the mapping $\di (0,\infty)\ni
t\longmapsto\frac{f(t)}{t}\quad\mbox{is nonincreasing};$
\medskip

\noindent $\di(f2)\qquad \lim_{t\ri 0^+}
\frac{f(t)}{t}=+\infty\quad\mbox{and}\;\;
\lim_{t\rightarrow\infty}\frac{f(t)}{t}=0.$
\medskip

We are concerned in this paper with  ground state solutions, that is, 
positive
solutions defined in the whole space and decaying to zero at
infinity. It is worth pointing out here that we do not assume any 
blow-up rate of decay on $g$ around the origin (as in Ghergu and R\u adulescu \cite{gr}) 
that may imply the property of compact support of B\'enilan,
Brezis and Crandall (see \cite{bbc}). 

There is a large number of works dealing with singular elliptic 
equations in bounded domains. In this sense,
we refer the reader to Callegari and Nachman \cite{cn}, C\ii rstea, Ghergu and R\u adulescu %%@
\cite{crg}, Coclite and Palmieri \cite{cp}, Crandall, Rabinowitz and Tartar \cite{crt}, D\'iaz, J.  %%@
Morel, and Oswald \cite{dmo}, Ghergu and R\u adulescu \cite {grcras}, Shi and Yao \cite{sy1}. The 
influence of the convection term has been 
emphasized in Ghergu and R\u adulescu \cite {gre,grjmaa} and Zhang \cite{zna}. Such singular 
boundary value problems
arise in the context of chemical heterogeneous catalysts and chemical
catalyst kinetics, in the
theory of heat conduction in electrically conducting materials,
singular minimal surfaces, as well
as in the study of non-Newtonian fluids or boundary layer phenomena for
viscous fluids. We
also point out that, due to the
meaning of the unknowns (concentrations, populations, etc.), only the
positive solutions are
relevant in most cases.

Concerning the ground state solutions for singular elliptic equations, 
we mention here the works of Lair and Shaker \cite{ls}, Sun and Liu \cite{sli}. In \cite{ls} it 
is considered the following singular boundary value problem 
\neweq{poo}
\left\{\begin{tabular}{ll}
$-\Delta u=p(x)g(u)$
\quad & $\mbox{\rm in}\ \RR^N,$\\
$u>0$ \quad & $\mbox{\rm in}\ \RR^N,$\\
$u(x)\ri 0$ \quad & $\mbox{\rm as}\ |x|\ri \infty\,,$
\end{tabular} \right.
\endeq
where $g:(0,\infty)\ri (0,\infty)$ is a smooth decreasing function (possibly, unbounded 
around the origin) such that $\int_0^1g(t)dt<\infty$. The results have been extended 
by C\ii rstea and R\u adulescu \cite{crjmaa} to the case where $g$ satisfies the weaker %%@
assumptions $\lim_{t\searrow 0}g(t)/t=+\infty$ and the mapping $t\longmapsto g(u)/(u+\beta)$
is decreasing, for some $\beta>0$.

It is proved in Lair and Shaker \cite{ls} (see also C\ii rstea and R\u adulescu \cite{crjmaa}) %%@
that a necessary condition in order to have a 
solution for \eq{poo} is
\neweq{ppw}
\di \int_1^\infty t\psi (t)dt<\infty,
\endeq
where $\psi(r)=\min_{|x|=r}p(x),$ $r\geq 0$. Note that condition 
\eq{ppw} is also necessary for our problem \eq{pp}, since any solution of \eq{pp} is a %%@
super-solution of 
\eq{poo}.
The sufficient condition for existence supplied in Lair and Shaker \cite{ls} is
\neweq{ppo}
\di \int_1^\infty t\phi (t)dt<\infty,
\endeq
where $\phi(r)=\max_{|x|=r}p(x),$ $r\geq 0$.
Hence, when $p$ is radially symmetric, then the problem \eq{poo}
has solutions if and only if $\int_1^\infty t p(t)dt<\infty$.

The main feature here is the presence of the convection
term $|\nabla u|^a$. In this sense we prove the following result.

\begin{theorem}\label{the1} Assume that $(f1)-(f2)$, $(g1)$ and 
\eq{ppo} are fulfilled. 
Then problem \eq{pp} has at least one solution.
 \end{theorem}

We  point out that the uniqueness of the solution to \eq{pp} is a delicate 
matter even in case of bounded domains (see, e.g., Ghergu and R\u adulescu \cite{gre,grjmaa}), 
due to the lack of an adequate comparison principle. We also notice that the growth decay
of the potential $p(x)$ described in our hypothesis \eq{ppo} implies that $p$ is in a certain
Kato class $K_{\rm loc}^N (\RR^N)$. This theory
was introduced by Aizenman and Simon in \cite{aize} 
to describe wide classes of functions arising in
Potential Theory.

\section{Proof of Theorem \ref{the1}}

The solution of problem \eq{pp} is obtained as a limit in
$C^{2,\gamma}_{{\rm loc}}(\RR^N)$ of a monotone sequence of
solutions associated to \eq{pp} in smooth bounded domains. A basic 
ingredient in our approach is the following auxiliary result.

\begin{lem}\label{p1} Let  $\Omega\subset\RR^N$ be a smooth bounded 
domain. Assume that
$f$ and $g$ satisfy $(f1)-(f2)$ and $(g1)$ respectively. Then the 
boundary value
problem
\begin{equation}\label{marginit}
\left\{\begin{tabular}{ll} $\di -\Delta
u=p(x)(g(u)+f(u)+|\nabla u|^a)$ \quad & $\mbox{\rm in}\ \Omega,$\\
$u>0$ \quad & $\mbox{\rm in}\ \Omega,$\\
$u=0$ \quad & $\mbox{\rm on}\ \partial\Omega,$\\
\end{tabular} \right.
\end{equation}
has a unique solution $u\in C^{2,\gamma}(\Omega)\cap
C(\overline\Omega).$
\end{lem}

\proof
The proof relies on sub and super-solution method. The assumptions on 
$f$ and $g$ imply that
 $m:=\inf_{t>0}\{g(t)+f(t)\}>0$. So, 
the unique solution $\underline u$ of the problem
\begin{equation}\label{zeta}
\left\{\begin{tabular}{ll}
$-\Delta u=m p(x)$ \quad & $\mbox{\rm in}\ \Omega,$\\
$\di u>0$ \quad & $\mbox{\rm in}\ \Omega,$\\
$\di u=0$ \quad & $\mbox{\rm on}\ \partial\Omega,$\\
\end{tabular} \right. \end{equation}
is a sub-solution of \eq{marginit}. 
The main point is to find a super-solution
$\overline u$ of problem (\ref{marginit}) such that $\underline u\leq
\overline u$ in $\Omega.$ Then, by classical results (see, e.g.,  
Gilbarg and Trudnger \cite{gt}) we deduce that the
problem (\ref{marginit}) has at least one solution.
\smallskip

Let $h:[0,\eta]\rightarrow[0,\infty)$ be such that
\begin{equation}\label{hh}
\left\{\begin{tabular}{ll}
$h''(t)=-g(h(t))\quad \mbox{ for all } 0<t<\eta,\,$\\
$h(0)=0,$\\
$h>0\quad\mbox{ in } (0,\eta].$
\end{tabular} \right.
\end{equation}
The existence of $h$ follows from the results in Agarwal and O'Regan 
\cite{ar}. 
Since $h$ is concave, there exists $h'(0+)\in (0,+\infty].$ 
Taking $\eta>0$ small enough, we can assume that $h'>0$ on $(0,\eta],$
that is, $h$ is increasing on $[0,\eta].$ Multiplying by $h'(t)$ in 
(\ref{hh})
and integrating on $[t,\eta],$ we obtain
\begin{equation}\label{ha}
\begin{aligned}
\di (h')^2(t)&\di =2\int_t^\eta
g(h(s))h'(s)ds+(h')^2(\eta)\\
&\di =2\int_{h(t)}^{h(\eta)} g(\tau)d\tau+(h')^2(\eta),\quad\mbox{ for 
all }\;0<t<\eta.
\end{aligned}
\end{equation}
Using the monotonicity of $g$ in (\ref{ha}) we get
$$\di (h')^2(t)\leq 2h(\eta)g(h(t))+(h')^2(\eta),\quad\mbox{ for all }
0<t<\eta.$$ Since $s^a\leq s^2+1,$ for all $s\geq 0,$
the last inequality yields
\begin{equation}\label{hprim}
\di (h')^a(t)\leq C g(h(t)),\quad\mbox{ for all }0<t<\eta,
\end{equation}
for some $C>0$. Let $\varphi_1$ be the normalized positive 
eigenfunction
corresponding to the first eigenvalue $\la_1$ of $-\Delta$ in 
$H^1_0(\Omega).$ 
We fix $c>0$ such that $c\|\varphi_1\|_{\infty}<\eta.$ 

By Hopf's maximum principle, there exist $\omega\subset\subset\Omega$ 
and
$\delta>0$ such that
\begin{equation}\label{omega}
\di |\nabla \varphi_1|>\delta \;\;\mbox{ in
}\;\Omega\setminus\omega \quad\mbox{ and 
}\;\;\varphi_1>\delta\quad\mbox{ in }\;\omega.
\end{equation}

Let $M>1$ be such that 
\neweq{bm2}
(Mc)^{1-a}\la_1 (h')^{1-a}(\eta)>3\max_{x\in\overline\Omega} p(x) 
\|\nabla \varphi_1\|^a_{\infty},
\endeq
\neweq{bm4}
Mc \la_1  h'(\eta)>3\max_{x\in\overline\Omega} p(x) 
g(h(c\min_{x\in\overline\omega}\varphi_1))
\endeq

\noindent and
\neweq{bm1}
\min\{ M(c\delta)^2, M^{1-a} C^{-1}(c\delta)^{2-a}\}> 
3\max_{x\in\overline\Omega} p(x),
\endeq
where $C$ is the constant from \eq{hprim}.
Since 
$$\di \lim_{t\ri 0^+} 
\Big((c\delta)^2g(h(t)-3\max_{x\in\overline\Omega}p(x) f(h(t)) \Big)=\infty,$$
we can assume that
\neweq{bm3}
\di (c\delta)^2 g(h(c\varphi_1))> 3 \max_{x\in\overline\Omega}p(x) 
f(h(c\varphi_1)) \quad\mbox{ in }\Omega\setminus\omega.
\endeq

Finally, from the assumption $(f2)$ on $f$ we have
$ \lim_{t\ri\infty}\frac{f(th(c\|\varphi_1\|_\infty)}{t}=0$,
so that we can choose $M>1$ large enough with the property
$$\di c\la_1 \inf_{x\in\overline\omega}\varphi_1 
h'(\eta)>\max_{x\in\overline\Omega}p(x) 
\frac{f(M h(c\|\varphi_1\|_\infty)}{M}.$$
The last inequality combined with the fact that $h'$ is decreasing 
yields
\neweq{bm5}
\di Mc\la_1\varphi_1 h'(c\varphi_1)\geq 3f(Mh(c\varphi_1))\quad\mbox{ 
in }\;\omega.
\endeq

We claim that $\overline u=Mh(c\varphi_1)$ is a super-solution of
the problem (\ref{marginit}) provided that $M$ satisfies
(\ref{bm2}), \eq{bm4} and (\ref{bm1}). We have 
$$ -\Delta
\overline u=Mc^2g(h(c\varphi_1))|\nabla\varphi_1|^2+
Mc\la_1\varphi_1h'(c\varphi_1)\quad\mbox{ in }\,\Omega.$$ 
From
(\ref{omega}), \eq{bm1} and the monotonicity of $g$ we obtain
\begin{equation}\label{ineg5}
\frac{1}{3}Mc^2g(h(c\varphi_1))|\nabla\varphi_1|^2\geq p(x)
g(h(c\varphi_1))\geq p(x)g(Mh(c\varphi_1))=p(x)g(\overline
u)\quad\mbox{ in }\;\Omega\setminus\omega.
\end{equation}
From (\ref{omega}), \eq{bm3} and our hypothesis $(f1)$, we obtain
\begin{equation}\label{ineg6}
\di \frac{1}{3}Mc^2g(h(c\varphi_1))|\nabla\varphi_1|^2\geq M p(x) 
f(h(c\varphi_1))
\geq p(x) f(M h(c\varphi_1))=p(x)f(\overline u)
\quad\mbox{ in }\;\Omega\setminus\omega.
\end{equation}
From \eq{hprim} and  \eq{bm1} we have
\begin{equation}\label{ineg6bis}
\di \frac{1}{3}Mc^2g(h(c\varphi_1))|\nabla\varphi_1|^2\geq p(x) 
(Mc h'(c\varphi_1) |\nabla \varphi_1|)^a=p(x)|\nabla \overline u|^a  
\quad\mbox{ in }\;\Omega\setminus\omega.
\end{equation}
Now, relations (\ref{ineg5})-(\ref{ineg6bis}) yield
\begin{equation}\label{ineg7}
-\Delta \overline v\geq
Mc^2g(h(c\varphi_1))|\nabla\varphi_1|^2\geq p(x)(g(\overline
u)+f(\overline u)+|\nabla \overline u|^a) \quad\mbox{ in
}\;\Omega\setminus\omega.
\end{equation}
Similarly, from \eq{bm2}-(\ref{bm5}) we deduce that
\begin{equation}{\label{ineg8}}
-\Delta \overline u\geq Mc\la_1\varphi_1h'(c\varphi_1)\geq
p(x)(g(\overline u)+f(\overline u)+|\nabla \overline u|^a)\quad\mbox{ 
in }\;
\omega.
\end{equation}
Using (\ref{ineg7}) and (\ref{ineg8}), it follows that $\overline u$ is 
a
super-solution of problem (\ref{marginit}). The maximum principle 
implies that 
$\underline u\leq \overline u$ in $\Omega$. Thus, the problem 
(\ref{marginit}) has
at least one classical solution. 
This concludes the proof of our Lemma.
\qed
\medskip

In what follows we apply Lemma \ref{p1} for $B_n:=\{x\in\RR^N; 
\,|x|<n\}.$
Hence, for all $n\geq 1$ there exists $u_n\in
C^{2,\gamma}(B_n)\cap C(\overline{B_n})$ such that
\begin{equation}\label{un}
\left\{\begin{tabular}{ll} $\di -\Delta u_n=p(x)(g(u_n)+f(u_n)+|\nabla
u_n|^a)$ \quad & $\mbox{\rm in}\ B_n,$\\
$u_n>0$ \quad & $\mbox{\rm in}\ B_n,$\\
$u_n=0$ \quad & $\mbox{\rm on}\ \partial B_n.$\\
\end{tabular} \right.
\end{equation}
We extend $u_n$ by zero outside of $B_n.$ We claim that 
$$\di u_n\leq u_{n+1}\quad \mbox{ in } B_n.$$
Assume by contradiction that the inequality $u_n\leq u_{n+1}$ does not
hold throughout $B_n$ and let 
$$\di \zeta(x)=\frac{u_n(x)}{u_{n+1}(x)},\quad x\in B_n.$$ Clearly
$\zeta=0$ on $\partial B_n$, so that $\zeta$ achieves its maximum in a 
point $x_0\in B_n$. At this point we have
$\nabla \zeta(x_0)=0$ and $\Delta \zeta(x_0)\leq 0$. This yields
$$\di -\mbox{div}(u_{n+1}^2\nabla 
\zeta)(x_0)=-\Big(\mbox{div}(u_{n+1}^2)\nabla \zeta+u_{n+1}^2\Delta \zeta  \Big)(x_0)\geq 0.$$
A straightforward computation shows that
$$\di -\mbox{div}(u_{n+1}^2\nabla \zeta)=-u_{n+1}\Delta u_n+u_n\Delta 
u_{n+1}.$$
Hence
$$\di \Big(-u_{n+1}\Delta u_n+u_n\Delta u_{n+1}\Big)(x_0)\geq 0.$$
The above relation produces
\neweq{cont1}
\di 
\left(\frac{g(u_n)+f(u_n)}{u_n}-\frac{g(u_{n+1})+f(u_{n+1})}{u_{n+1}}\right)(x_0)+
\left(\frac{|\nabla u_{n}|^a}{u_n}-\frac{|\nabla 
u_{n+1}|^a}{u_{n+1}}\right)(x_0)\geq 0.
\endeq

Since $t\longmapsto \frac{g(t)+f(t)}{t}$ is decreasing on $(0,\infty)$ 
and $u_n(x_0)>u_{n+1}(x_0)$, from \eq{cont1} 
we obtain

\neweq{cont2}
\di \left(\frac{|\nabla u_{n}|^a}{u_n}-\frac{|\nabla 
u_{n+1}|^a}{u_{n+1}}\right)(x_0)> 0.
\endeq 
On the other hand, $\nabla \zeta(x_0)=0$ implies
$$\di u_{n+1}(x_0)\nabla u_n(x_0)=u_n(x_0)\nabla u_{n+1}(x_0).$$
Furthermore, relation \eq{cont2} leads us to
$$\di u_n^{a-1}(x_0)-u^{a-1}_{n+1}(x_0)>0,$$
which is a contradiction since $0<a<1$.
Hence $u_n\leq u_{n+1}$ in $B_n$ which means that 
$$\di 0\leq u_1\leq \dots\leq u_n\leq u_{n+1}\leq \dots \quad\mbox{ in 
}\;\RR^N.$$
The main point is to find an upper bound for the sequence $(u_n)_{n\geq 
1}$. 
This is provided in the following result.

\begin{lem}\label{ll1}
 The inequality problem
\neweq{ppq}
\left\{\begin{tabular}{ll}
$-\Delta v\geq p(x)(g(v)+f(v)+ |\nabla v|^a)$
\quad & $\mbox{\rm in}\ \RR^N,$\\
$v>0$ \quad & $\mbox{\rm in}\ \RR^N,$\\
$v(x)\ri 0$ \quad & $\mbox{\rm as}\ |x|\ri \infty$
\end{tabular} \right.
\endeq
has at least one solution in $C^2(\RR^N)$.
\end{lem}

\proof Set 
$$\di \Phi(r)=r^{1-N}\int_0^rt^{N-1}\phi(t)dt, \quad\mbox{ for all }\; 
r>0.$$ 
Using the assumption (\ref{ppo}) and L'H\^opital's rule, we get
$\lim_{r\ri\infty}\Phi(r)=\lim_{r\searrow 0}\Phi(r)=0.$ Thus,
$\Phi$ is bounded on $(0,\infty)$ and it can be extended in the origin 
by
taking $\Phi(0)=0.$ On the other hand, we have
$$\begin{aligned}
\di \int_0^r\Phi(t)dt&\di =-\frac{1}{N-2}\int_0^r
\frac{d}{dt}\left(t^{2-N}\right)\int_0^t s^{N-1}\phi(s)dsdt\\
&\di =-\frac{1}{N-2}\left[ r^{N-2}
\int_0^r s^{N-1}\phi(s)ds-\int_0^r t\phi(t)dt\right]\\
&\di =\frac{1}{N-2}\frac{\di-\int_0^r s^{N-1}\phi(s)ds+
r^{N-2}\int_0^r t\phi(t)dt}{r^{N-2}}.\\
\end{aligned}$$

By L'H\^opital's rule we obtain
\neweq{aalpha}\di \lim_{r\ri\infty}\frac{\di-\int_0^r s^{N-1}\phi(s)ds+
r^{N-2}\int_0^r t\phi(t)dt}{r^{N-2}}=\int_0^\infty r\phi(r)dr.\endeq 
The
last two relations imply
$$
\di\int_0^\infty\Phi(r)dr=\di\lim_{r\ri\infty}\int_0^r
\Phi(t)dt=\frac{1}{N-2}\int_0^\infty r\phi(r)dr<\infty.
$$

Let $k>2$ be such that
\begin{equation}\label{kk}
\di k^{1-a}\geq 2\max_{r\geq 0}\Phi^a(r).
\end{equation}
In view of \eq{aalpha} we can define 
$$\di \xi(x)=k\int_{|x|}^\infty\Phi(t)dt,\quad\mbox{ for all 
}\;x\in\RR^N.$$ 
Then $\xi$ satisfies
$$\left\{\begin{tabular}{ll}
$-\Delta \xi=k\phi(|x|)$ \quad & $\mbox{\rm in}\ \RR^N,$\\
$\xi>0$ \quad & $\mbox{\rm in}\ \RR^N,$\\
$\xi(x)\ri0$ \quad & $\mbox{\rm as}\ |x|\ri\infty.$
\end{tabular}\right.$$
Since the mapping $ [0,\infty)\ni t\longmapsto
\int_0^t\frac{1}{g(s)+1}ds\in[0,\infty)$ is bijective, we can
implicitly define $w:\RR^N\ri(0,\infty)$ by
$$\di \int_0^{w(x)}\frac{1}{g(t)+1}dt=\xi(x),\quad \mbox{
 for all }\;x\in\RR^N.$$ 
It is easy to see that $w\in C^2(\RR^N)$ and $w(x)\ri 0$ as $|x|\ri 
\infty.$ 
Furthermore, we have
\neweq{ctrd}
\di |\nabla w|=|\nabla \xi|(g(w)+1)=k\Phi(|x|)(g(w)+1)\quad \mbox { in 
}\,\RR^N
\endeq
and
$$\begin{aligned}
\di -\Delta w&\di =-(g(w)+1)\Delta \xi-g'(w)(g(w)+1)|\nabla \xi|^2
 \geq k\phi(|x|)(g(w)+1)\\
&\di \geq \phi(|x|)(g(w)+1)+\frac{1}{2}k\phi(|x|)(g(w)+1).\\
\end{aligned}$$
By (\ref{kk}) and \eq{ctrd} we deduce
$$ 
\di \frac{k}{2}\phi(|x|)(g(w)+1)\geq
\phi(|x|)k^a(g(w)+1)^a \Phi^a(|x|)\geq p(x) |\nabla w|^a \quad\mbox{
in }\;\RR^N.$$
 
Hence
\neweq{mlt}
\left\{\begin{tabular}{ll}
$-\Delta w\geq p(x)(g(w)+1+|\nabla w|^a)$ \quad & $\mbox{\rm in}\ 
\RR^N,$\\
$w>0$ \quad & $\mbox{\rm in}\ \RR^N,$\\
$w(x)\ri 0$ \quad & $\mbox{\rm as}\ |x|\ri\infty.$\\
\end{tabular}\right.
\endeq
Using the assumption $(f1)$, we can find $M>1$ large enough such that 
$M>f(Mw)$ in $\RR^N$. Multiplying by $M$ in
\eq{mlt} we deduce that $v:=Mw$ satisfies \eq{ppq} and the proof of 
Lemma \ref{ll1} is now complete. 
\qed
\medskip

{\bf Proof of Theorem \ref{the1} concluded.}
With the same proof as above we deduce that $u_n\leq v$ in $B_n$, for 
all $n\geq 1$. This implies 
$$\di 0\leq u_1\leq \dots\leq u_n\leq v \quad\mbox{ in }\;\RR^N.$$

Thus,
there exists $u(x)=\lim_{n\ri\infty}u_n(x),$ for all $x\in\RR^N$ and 
$u_n\leq
u\leq v$ in $\RR^N.$ Since $v(x)\ri 0$ as $|x|\ri \infty,$ we
deduce that $u(x)\ri 0$ as $|x|\ri \infty.$ 
A standard bootstrap argument (see Gilbarg and Trudinger \cite{gt}) 
implies that
 $\di u_n\ri u$ in $C^{2,\gamma}_{{\rm
loc}}(\RR^N)$ and that $u$ is a solution of problem \eq{pp}.

This finishes the proof of Theorem \ref{the1}.
\qed


\begin{thebibliography}{99}  {\footnotesize

\bibitem{ar} R. Agarwal and D. O'Regan, Existence theory for single
and multiple solutions to singular positone boundary value problems,
{\it J. Differential Equations} {\bf 175} (2001), 393-414.

\bibitem{aize} M. Aizenman and B. Simon, Brownian motion and Harnack inequality for 
Schr\"odinger operators, {\it Comm. Pure Appl. Math.} {\bf 35} (1982), 209-273.

\bibitem{ap} N. E. Alaa and M. Pierre, Weak solutions of
some quasilinear elliptic equations with data measures, {\it SIAM
J. Math. Anal.} {\bf 24} (1993), 23-35.

\bibitem{arcoya} D. Arcoya, Positive solutions for semilinear Dirichlet problems in an annulus,
{\it J. Differential Equations} {\bf 94} (1991),
217-227.

\bibitem{bbc} P. B\'enilan, H. Brezis, and M. Crandall,
A semilinear equation in $L^1(\RR^N)$, {\it Ann. Scuola Norm. Sup.
Pisa} {\bf 4} (1975), 523-555.

\bibitem{cn} A. Callegari and A. Nachman, A nonlinear
singular boundary value
problem in the theory of pseudoplastic  fluids, {\it
SIAM J. Appl. Math.} {\bf 38} (1980), 275-281.

\bibitem {crg} F. C\ii rstea, M. Ghergu, and V. R\u adulescu,
Combined effects of asymptotically linear and singular
nonlinearities in bifurcation problems of Lane-Emden-Fowler type,
{\it  J. Math. Pures Appliqu\'ees} {\bf 84}
(2005), 493-508.

\bibitem {crjmaa} F. C\ii rstea and V. R\u adulescu,
Existence and uniqueness of positive solutions to a semilinear elliptic problem in $\RR^N$,
{\it  J. Math. Anal. Appl.} {\bf 229} (1999), 417-425.

\bibitem{cp} M. M. Coclite and G. Palmieri, On a
singular nonlinear Dirichlet problem, {\it Comm. Partial
Differential Equations} {\bf 14} (1989), 1315-1327.

\bibitem{crt} M. G. Crandall, P. H. Rabinowitz and L. Tartar, On a
Dirichlet problem
with a singular nonlinearity, {\it Comm. Partial Differential
Equations} {\bf 2} (1977), 193-222.

\bibitem{dmo} J. I. D\'iaz, J. M. Morel, and L. Oswald,
An elliptic equation with singular nonlinearity, {\it Comm. Partial
Differential
Equations} {\bf 12} (1987), 1333-1344.

\bibitem{fowler} R. H. Fowler, The solution of Emden's and similar differential equations, 
{\it Monthly Notices Roy. Astro. Soc.} {\bf 91} (1930), 63-91.

\bibitem{gr} M. Ghergu and V. R\u adulescu, Sublinear
singular elliptic problems with two parameters,
{\it J.~Differential Equations} {\bf 195} (2003), 520-536.

\bibitem{grcras} M. Ghergu and V. R\u adulescu, Bifurcation for a class 
of singular elliptic problems with quadratic convection term, {\it C. R. Acad. Sci. 
Paris, Ser. I} {\bf  338} (2004), 831-836.

\bibitem{gre} M. Ghergu and V. R\u adulescu, Multiparameter
bifurcation and asymptotics for the singular Lane-Emden-Fowler
equation with convection term, {\it Proc. Royal Soc. Edinburgh
Sect. A} {\bf 135} (2005), 61-84.

\bibitem {grjmaa} M. Ghergu and V. R\u adulescu, On a class of
sublinear
singular elliptic problems with convection term, {\it J.~Math.
Anal. Appl.} {\bf 311} (2005), 635-646.

\bibitem{gt} D. Gilbarg and N. S. Trudinger, {\it Elliptic Partial
Differential Equations of Second Order}, 2nd ed., Springer-Verlag,
Berlin Heidelberg New York, 1983.

\bibitem{hale} J. K. Hale and J. Mawhin, Coincidence degree and periodic solutions of 
neutral equations, {\it J. Differential Equations} {\bf 15} (1974), 295-307.

\bibitem{ls} A. V. Lair and A. W. Shaker, Classical and weak solutions
of a singular semilinear elliptic problem, {\it J.~Math.
Anal. Appl.} {\bf 211} (1997), 371-385.

\bibitem{mw} J. Mawhin and M. Willem, {\it Critical Point Theory and Hamiltonian 
Systems}, Springer, New York, 1989.

\bibitem{nato} Y. Nato and S. Tanaka, On the existence of multiple solutions of the boundary 
value problems for nonlinear second order
differential equations, {\it Nonlinear Anal.} {\bf 56} (2004), 919-935.

\bibitem {sy1} J. Shi and M. Yao, On a singular
nonlinear semilinear elliptic problem, {\it Proc. Royal Soc.
Edinburgh Sect. A} {\bf 128} (1998), 1389-1401.

\bibitem{sli} Y. Sun and S. Li, Structure of ground state solutions 
of singular semilinear elliptic equation, {\it Nonlinear Analysis} {\bf 
55} (2003), 399-417. 

\bibitem{wang} H. Wang, On the existence of positive solutions for semilinear elliptic 
equations in annulus, {\it J. Differential Equations} {\bf 109} (1994), 1-7.

\bibitem{wong} J. S. W. Wong, On the generalized Emden-Fowler equation, {\it SIAM Rev.} 
{\bf 17} (1975), 339-360.

\bibitem{zna} Z. Zhang, Nonexistence of positive classical solutions of
a singular nonlinear
Dirichlet problem with a convection term, {\it Nonlinear Anal., T.M.A.}
{\bf 27} (1996), 957-961.






 }


\end{thebibliography}
\end{document}